\documentclass[9pt]{amsart}
\usepackage{amsmath}
\usepackage{amsfonts}
\usepackage{graphicx}
\usepackage{amscd}
\usepackage{amsmath,amsfonts,amssymb,amsthm,mathrsfs,xypic}
\newtheorem{thm}{Theorem}[section]

\newtheorem{cor}[thm]{Corollary}
\newtheorem{lem}[thm]{Lemma}
\newtheorem{exm}[thm]{Example}

\newtheorem{defn}[thm]{Definition}
\numberwithin{equation}{section}

\begin{document}

\title[Gorensteinness, homological invariants and Gorenstein derived categories]
{Gorensteinness, homological invariants and Gorenstein derived categories}
\author[Nan Gao, ]{Nan  Gao}
\thanks{ e-mail: gaonanjane$\symbol{64}$gmail.com, nangao$\symbol{64}$shu.edu.cn}
\thanks{Supported by the National Natural Science Foundation of China
(Grant No. 11101259).}

\maketitle
\begin{center}
Department of Mathematics, \ \ Shanghai  University\\
Shanghai 200444, P. R. China
\end{center}

\vskip 5pt

\begin{abstract} Relations between Gorenstein derived categories, Gorenstein defect categories and Gorenstein stable categories are established.
Using these, the Gorensteinness of an algebra $A$ and invariants with respect to recollements of the bounded Gorenstein derived category $D^{b}_{gp}(A\mbox{-}{\rm mod})$ of $A$ are investigated. Specifically, the Gorensteinness of  $A$ is characterized in three ways: the existence of Auslander-Reiten triangles in  $D^{b}_{gp}(A\mbox{-}{\rm mod})$;  recollements of $D^{b}_{gp}(A\mbox{-}{\rm mod})$;  and also Gorenstein derived equivalences. It is shown that  the finiteness of Cohen-Macaulay type  and of finitistic dimension are invariant with respect to the recollements of $D^{b}_{gp}(A\mbox{-}{\rm mod})$.

\vskip 10pt


\vskip 10pt

\end{abstract}

\vskip 20pt

\section{\bf Introduction}

\vskip 10pt

The Gorensteinness of an algebra is of interest in the representation theory of
algebras, in Gorenstein homological algebra, and in the theory of singularity categories
(see e.g. [Ha2], [CV], [Be2], [BJO], [K], [LZ]). How to characterize Gorenstein property is a basic problem.

\vskip 10pt

An algebra has many invariants, for example, the finiteness of global dimension, the finiteness of finitistic dimension, and the finiteness of Cohen-Macaulay type (see e.g. [Wi], [Ha4], [Be4]). How to describe and compare these homological invariants is a major topic of interest.

\vskip 10pt

One can approach these questions by derived (and related) categories, as well as comparisons of derived categories.
Derived categories, introduced by Grothendieck and Verdier([Ver]), have been playing an increasingly important role in various areas of mathematics, including
representation theory, algebraic geometry, and mathematical physics. For example, Happel ([Ha3])  has shown  a finite dimensional algebra $A$ has finite global dimension if and only if the bounded derived category $D^{b}(A\mbox{-}{\rm mod})$ has Auslander-Reiten triangles.
There are two ways to compare derived categories. One way is by derived equivalences. Happel ([Ha1]) has shown the finiteness of global dimension of an algebra is invariant under derived equivalences.
Another way is by  recollements, which have been introduced by Beilinson, Bernstein and Deligne ([BBD]). A recollement of a derived category by another two derived categories is a diagram of six functors between these categories, generalising Grothendieck's six functors. Suppose that $A$, $B$ and $C$ are three finite dimensional algebras over a field. If  $D^{b}(A\mbox{-}{\rm mod})$ admits a recollement with respect to the bounded derived categories $D^{b}(B\mbox{-}{\rm mod})$ and $D^{b}(C\mbox{-}{\rm mod})$ of $B$ and $C$,   Wiedemann ([Wi]) has shown that  $A$ has finite global dimension if and only if $B$ and $C$ have so, and Happel ([Ha4]) has shown that  $A$ has finite finitistic dimension if and only if $B$ and $C$ have so.

\vskip 10pt

For Gorenstein homological algebra we refer to [AM, Buc, CFH,  EJ, Hol].
Nan Gao and Pu Zhang defined the corresponding version of the derived
category in Gorenstein homological algebra. They introduced in [GZ]
the notions of Gorenstein derived category and Gorenstein derived equivalence which are needed
in this context. Following [GZ], the  bounded Gorenstein derived category
$D^b_{gp}(A\mbox{-}{\rm mod})$ of an algebra $A$ is
defined as the Verdier quotient of the bounded homotopy category
$K^b(A\mbox{-}{\rm mod})$ with respect to the triangulated
subcategory $K^b_{gpac}(A\mbox{-}{\rm mod})$ of
$A\mbox{-}\mathcal Gproj$-acyclic complexes. Later, Gao ([G2]) described  Auslander-Reiten triangles in the bounded Gorenstein derived categories for Gorenstein algebras of finite CM-type. In [G3] a sufficient and necessary  criterion
has been given for the existence of  recollements of Gorenstein derived categories.
Based on these work, two questions arise.

\vskip 5pt

$(1)$ \  Can we characterize the Gorensteinness
of an algebra in terms of  the corresponding Gorenstein derived category?

\vskip 5pt

$(2)$ \  Which invariants can be compared along recollements of Gorenstein derived categories?

\vskip 10pt

In this article we will provide answers to these questions. Our answer to question $(1)$
are the combination of  Corollary 3.2 and 3.3, and Theorem 3.4 and 3.9(2). We state them as Theorem $A$.

\vskip 10pt

{\bf Theorem A} \ Let $A$ be an artin algebra. Consider the following statements:

\vskip 5pt

$(1)$ \ $A$ is Gorenstein;

\vskip 5pt

$(2)$ \  $A$ is virtually Gorenstein and $D^{b}_{gp}(A\mbox{-}{\rm mod})$ has Auslander-Reiten triangles;

\vskip 5pt

$(3)$ \  $A$ is virtually Gorenstein, and there exist Gorensein algebras $B$ and $C$ and a recollement 
$$
D^{b}_{gp}(B\mbox{-}{\rm mod}) \begin{smallmatrix}
  \underleftarrow{ \ \ \  i^{*} \ \ \  } \\
  \underrightarrow{ \ \ \ i_{*} \ \ \ } \\
  \overleftarrow{ \ \ \ i^{!} \ \ \ }
\end{smallmatrix}
D^{b}_{gp}(A\mbox{-}{\rm mod})
\begin{smallmatrix}
  \underleftarrow{ \ \ \ j_{!} \ \ \ } \\
  \underrightarrow{\ \ \ j^{*} \ \ \ } \\
  \overleftarrow{ \ \ \ j_{*} \ \ \ }
\end{smallmatrix}
D^{b}_{gp}(C\mbox{-}{\rm mod});
$$

$(3^{'})$ \  $A$ is virtually Gorenstein, and for arbitrary virtually  Gorensein algebras $B$ and $C$, if there exists the following recollement 
$$
D^{b}_{gp}(B\mbox{-}{\rm mod}) \begin{smallmatrix}
  \underleftarrow{ \ \ \  i^{*} \ \ \  } \\
  \underrightarrow{ \ \ \ i_{*} \ \ \ } \\
  \overleftarrow{ \ \ \ i^{!} \ \ \ }
\end{smallmatrix}
D^{b}_{gp}(A\mbox{-}{\rm mod})
\begin{smallmatrix}
  \underleftarrow{ \ \ \ j_{!} \ \ \ } \\
  \underrightarrow{\ \ \ j^{*} \ \ \ } \\
  \overleftarrow{ \ \ \ j_{*} \ \ \ }
\end{smallmatrix}
D^{b}_{gp}(C\mbox{-}{\rm mod}),
$$
then $B$ and $C$ are Gorenstein;

\vskip 5pt

$(4)$ \  There is a triangle-equivalence $D^{b}_{gp}(A\mbox{-}{\rm{mod}})\cong K^{b}(A\mbox{-}\mathcal{G}proj)$;

\vskip 5pt

$(5)$ \  $D^{b}_{gp}(A\mbox{-}{\rm mod})$ has Auslander-Reiten triangles;

\vskip 5pt

$(6)$ \  $A$ is  Gorenstein derived equivalent to an algebra $B$, which is  Gorenstein.

\vskip 5pt

\noindent We have the following relations between these statements:

\vskip 5pt

$(i) \ (1)\Longleftrightarrow (2)\Longleftrightarrow (3)\Longleftrightarrow(3^{'})$

\vskip 5pt

$(ii) \ (1)\Longleftrightarrow (6)$

\vskip 5pt

$(iii)$ \ If $A\mbox{-}\mathcal{G}proj$ is contravariantly
finite in $A\mbox{-}{\rm mod}$, then $(1)\Longleftrightarrow (4)$.

\vskip 5pt

$(iv)$ \ If $A$ is of finite CM-type, then $(1)\Longleftrightarrow (5)$.

\vskip 10pt

Our answer to question $(2)$ are the combination of Theorem 3.9(1) and 3.12. We state them as Theorem $B$.

\vskip 10pt

{\bf Theorem B} \ Let $A$, $B$ and $C$ be artin algebras. Assume that the bounded Gorenstein derived category $D^{b}_{gp}(A\mbox{-}{\rm mod})$ admits a
recollement with respect to $D^{b}_{gp}(B\mbox{-}{\rm mod})$ and $D^{b}_{gp}(C\mbox{-}{\rm mod})$. The following hold true:

\vskip 5pt

$(1)$ \ If $A$, $B$ and $C$ are virtually Gorenstein, then $A$ is of finite CM-type if and only if $B$ and $C$ are so;

\vskip 5pt

$(2)$ \ If $A$, $B$ and $C$ are of finite CM-type, then ${\rm fd}A<\infty$ if and only if
\  ${\rm fd}B<\infty$ and ${\rm fd}C<\infty$.

\vskip 10pt

Let us end this introduction by mentioning that in private communication with Javad Asadollahi, he points out that he and his collaborators  also
have proofs for Corollary 3.2 and Theorem 3.9 in this paper. Their proofs are obtained independently, and also
are different from the proofs given in this paper.  The author would like to thank him for letting us know their proofs.

\vskip 20pt

\section{\bf Preliminaries}

\vskip 10pt

In this section we fix notation and recall the main concepts to be used.

\vskip 10pt

Let $A$ be an  artin algebra. Denote by $A\mbox{-}{\rm Mod}(\rm
{resp}. \ A\mbox{-}{\rm mod})$ the category of  left
$A\mbox{-}$modules (resp. the category of finitely-generated left
$A\mbox{-}$modules), and $A\mbox{-}{\rm Proj}(\rm {resp.} $ $\
A\mbox{-}{\rm proj})$ the full subcategory of projective
$A\mbox{-}$modules (resp. the full subcategory of finitely-generated
projective $A\mbox{-}$modules). An $A$-module $M$ is said to be
Gorenstein-projective in $A$-Mod(resp. $A$-mod), if there is an
exact sequence $P^{\bullet}=\cdots \longrightarrow
P^{-1}\longrightarrow P^{0} \stackrel{d^0}{\longrightarrow}
P^{1}\longrightarrow P^{2}\longrightarrow \cdots$ in $A\mbox{-}{\rm
Proj}(\rm {resp.} \ A\mbox{-}{\rm proj})$ with
$\rm{Hom}_A(P^{\bullet}, Q)$ exact for any $A\mbox{-}$module $Q$ in
$A\mbox{-}{\rm Proj}(\rm {resp.} \ A\mbox{-} proj)$, such that
$M\cong \operatorname{ker}d^0$ (see [EJ]). Denote by
$A\mbox{-}\mathcal{GP}({\rm resp.} \
A\mbox{-}\mathcal{G}proj)$ the full subcategory of
Gorenstein-projective modules in $A\mbox{-}{\rm Mod}(\rm {resp}. \
A\mbox{-}{\rm mod})$,  and similarly denote by
$A\mbox{-}\mathcal{GI}$ the full subcategory of
Gorenstein-injective modules in $A\mbox{-}{\rm Mod}$.

\vskip 10pt

A proper Gorenstein-projective resolution of $A$-module $M$ in $A\mbox{-}{\rm mod}$ is an
exact sequence $E^{\bullet}=\cdots\longrightarrow G_1\longrightarrow
G_0\longrightarrow M\longrightarrow 0$ such that all $G_i\in
A\mbox{-}\mathcal Gproj$, and that
$\operatorname{Hom}_A(G, E^\bullet)$ stays exact for each
$G\in A\mbox{-}\mathcal Gproj$. The second requirement
guarantees the uniqueness of such a resolution in the homotopy
category (the Comparison Theorem; see [EJ], p.169).

\vskip 10pt

 The Gorenstein-projective dimension $\mathcal{G} p{\rm
dim}M$ of $M$ in $A\mbox{-}{\rm mod}$ is defined to be the smallest integer $n\geq 0$ such
that there is an exact sequence $0\longrightarrow
G_n\longrightarrow\cdots \longrightarrow G_1\longrightarrow
G_0\longrightarrow M\longrightarrow 0$ with all $G_i\in
A\mbox{-}\mathcal Gproj$, if it exists; and $\mathcal{G} p{\rm dim}M = \infty$ if there is no such exact sequence
of finite length. For an $A\mbox{-}$module $_{A}X$ we denote by ${\rm proj.dim}_{A}X$ the
projective dimension of $X$. Clearly $\mathcal{G}p{\rm dim}M \le {\rm
proj.dim}_{A}M$. Denote by ${\rm fGd}(A)={\rm sup}\{\mathcal{G}p{\rm dim}_{A}X|$ $\mathcal{G}p{\rm dim}_{A}X<\infty\}$.

\vskip 10pt

A complex $C^\bullet$ of finitely-generated $A\mbox{-}$modules is  $A\mbox{-}\mathcal Gproj$-acyclic,
if $\operatorname{Hom}_{ A}(G, C^\bullet)$ is acyclic for
any $G\in A\mbox{-}\mathcal Gproj$. It is also called  proper
exact for example in [AM]. A chain map \ $f^\bullet: X^\bullet
\longrightarrow Y^\bullet$ is
 an $A\mbox{-}\mathcal Gproj$-quasi-isomorphism, if \
 $\operatorname{Hom}_{ A}(G, f^\bullet)$ is a
 quasi-isomorphism for any $G\in A\mbox{-}\mathcal Gproj$, i.e., there
 are isomorphisms of abelian groups
 ${\rm H}^n\operatorname{Hom}_{ A}(G, f^\bullet): \
 {\rm H}^n\operatorname{Hom}_{ A}(G, X^\bullet)\cong
 {\rm H}^n\operatorname{Hom}_{A}(G, Y^\bullet), \ \forall \
 n\in\Bbb Z, \ \forall \ G\in A\mbox{-}\mathcal Gproj$.

 \vskip 10pt

 Denote by $K^{-}(A\mbox{-}\mathcal{G}proj)$  the upper bounded homotopy category of
$A\mbox{-}\mathcal{G}proj$, and by $K^{-, gpb}(A\mbox{-}\mathcal{G}proj)$ the full subcategory
\begin{align*}
K^{-, gpb}(A\mbox{-}\mathcal{G}proj):&=\{G^{\bullet}\in K^{-}(A\mbox{-}\mathcal{G}proj)\mid \ \exists \ a \ positive \ integer \ N \
 such \ that \\
&H^{-n}{\rm Hom}_{A}(E, G^{\bullet})=0, \ \forall n> N, \  \forall E\in A\mbox{-}\mathcal{G}proj\}.
\end{align*}
Note from [GZ, Theorem 3.6] that if $A\mbox{-}\mathcal Gproj$ is contravariantly finite in $A\mbox{-}{\rm mod}$, then there is a triangle-equivalence $D^{b}_{gp}(A\mbox{-}{\rm mod})\cong K^{-, gpb}(A\mbox{-}\mathcal{G}proj)$.

\vskip 10pt

We say that two artin algebras $A$ and $B$ are Gorenstein
derived equivalent, if there is a triangle-equivalence $D^b_{gp}(A\mbox{-}{\rm mod})\cong D^b_{gp}(B\mbox{-}{\rm mod})$.

\vskip 10pt

Recall from [BH,\ Be4] that an artin algebra $A$ is of finite Cohen-Macaulay type (simply, CM-type) if there are only finitely many isomorphism classes of finitely-generated indecomposable Gorenstein-projective $A$-modules. Suppose $A$ is an artin algebra of finite CM-type, and $G_1, \cdots, G_n$ are
all the pairwise non-isomorphic indecomposable finitely-generated
Gorenstein-projective $A$-modules, and $G = \bigoplus\limits_{1\le
i\le n}G_i$. Set $\mathcal{G}p(A):={\rm End}_{A}(G)^{\rm op}$, which
we call the relative Auslander algebra of $A$. It is clear that $G$
is an $A\mbox{-}\mathcal{G}p(A)$-bimodule and $\mathcal{G}p(A)$ is
an artin algebra ([ARS, p.27]). Denote by ${\rm{mod}}\mbox{-}\mathcal{G}p(A)$ the category of left $\mathcal{G}p(A)$-modules.
Recall from [BR, \ Be3] that an artin algebra $A$ is called virtually Gorenstein if
$A\mbox{-}\mathcal{GP}^{\perp}=^{\perp}A\mbox{-}\mathcal{GI}$.

\vskip 10pt

Let $A$ be an  artin algebra. Denote by $A\mbox{-}\mathcal{G}proj(M, N)$ the subgroup of ${\rm Hom}_A(M, N)$
of $A$-maps from $M$ to $N$ which factors through the finitely-generated Gorenstein-projective modules, and
$A\mbox{-}{\rm{mod}}/A\mbox{-}\mathcal{G}proj$ the quotient
category of $A\mbox{-}{\rm{mod}}$ modulo $A\mbox{-}\mathcal{G}proj$,
i.e., the objects are same as those of $A\mbox{-}{\rm{mod}}$, and the morphism space
from $M$ to $N$  is the quotient group ${\rm Hom}_A(M, N)/A\mbox{-}\mathcal{G}proj(M,
N)$. In the following, we call $A\mbox{-}{\rm{mod}}/A\mbox{-}\mathcal{G}proj$ the Gorenstein stable category of $A$.

\vskip 10pt

Following [Be1, Definition 3.1], the stabilization of $A\mbox{-}{\rm{mod}}/A\mbox{-}\mathcal{G}proj$
is a pair $(S, S(A\mbox{-}{\rm{mod}}$ $/A\mbox{-}\mathcal{G}proj))$, where $S(A\mbox{-}{\rm{mod}}/A\mbox{-}\mathcal{G}proj)$
is a triangulated category and $S: A\mbox{-}{\rm{mod}}/A\mbox{-}\mathcal{G}proj$ $\to S(A\mbox{-}{\rm{mod}}/A\mbox{-}\mathcal{G}proj)$ is an exact functor,
such that for any exact functor $F: A\mbox{-}{\rm{mod}}/A\mbox{-}\mathcal{G}$ $proj\to \mathcal{C}$ to a triangulated category $\mathcal{C}$, there exists a unique triangle functor
 $F^{*}: S(A\mbox{-}{\rm{mod}}/A\mbox{-}\mathcal{G}proj)\to \mathcal{C}$ such that $F^{*}S=F$. For the construction of $S(A\mbox{-}{\rm{mod}}/A\mbox{-}\mathcal{G}proj)$,
 see ([He], [Be1]).

\vskip 20pt

\section{\bf  On Gorenstein derived categories}

\vskip 10pt

In this section, we characterize the Gorensteinness of an algebra $A$ in three ways: the existence of Auslander-Reiten triangles in the bounded Gorenstein derived category  $D^{b}_{gp}(A\mbox{-}{\rm mod})$ of $A$;  recollements of $D^{b}_{gp}(A\mbox{-}{\rm mod})$,  and also Gorenstein derived equivalences.

\vskip 10pt

Let $A$ be an artin algebra and $D_{sg}(A)$ the singularity category of $A$.
Buchweitz's Theorem ([Buc]) has shown that there is an embedding $F:
\underline{A\mbox{-}\mathcal{G}proj}\to D_{sg}(A)$ given by
$F(G)=G$, where the second $G$ is the corresponding stalk complex at
degree $0$, and that if $A$ is Gorenstein, then $F$ is a
triangle-equivalence. The converse is also true (cf. Theorem 6.9 in [Be1]).  In general, to measure how far a ring is from being Gorenstein,
Bergh, J\o rgensen and Oppermann ([BJO]) defined the Gorenstein defect category
$D^{b}_{defect}(A):=D_{sg}(A)/ {\rm Im}F$,  and also they have shown that
$A$ is Gorenstein if and only if $D^{b}_{defect}(A)=0$.

\vskip 10pt

To recognize Gorenstein rings via Gorenstein derived categories, we compare a Gorenstein defect category with a Gorenstein derived category, and try to construct  a precise relation.  We start with the following lemma.

\vskip 10pt

\begin{lem} \  Let $A$ be an artin algebra such that $A\mbox{-}\mathcal{G}proj$ is contravariantly
finite in $A\mbox{-}{\rm mod}$. Then there is a triangle-equivalence  $$D^{b}_{defect}(A)\cong D^{b}_{gp}(A\mbox{-}{\rm mod})/K^{b}(A\mbox{-}\mathcal{G}proj).$$
\end{lem}

\vskip 5pt

\noindent{\bf Proof.}\quad By [KZ, Final Remark] we have a triangle-equivalence
$$D^{b}_{defect}(A)\cong K^{-gpb}(A\mbox{-}\mathcal{G}proj)/K^{b}(A\mbox{-}\mathcal{G}proj).$$
By [GZ, Theorem 3.6] we get that
$$D^{b}_{gp}(A\mbox{-}{\rm mod})\cong K^{-gpb}(A\mbox{-}\mathcal{G}proj).$$
This completes the proof.

\vskip 10pt

Now we  test the Gorensteinness of $A$ by the structures of $D^{b}_{gp}(A\mbox{-}{\rm mod})$.

\vskip 10pt

\begin{cor} \ Let $A$ be an artin algebra such that $A\mbox{-}\mathcal{G}proj$ is contravariantly
finite. Then $A$ is  Gorenstein if and only if there is a triangle-equivalence $D^{b}_{gp}(A\mbox{-}{\rm{mod}})\cong K^{b}(A\mbox{-}\mathcal{G}proj)$.
\end{cor}

\vskip 5pt

\noindent{\bf Proof.}\quad By [GZ, Corollary 3.8] we only need to verify the sufficiency. Since $D^{b}_{gp}(A\mbox{-}{\rm{mod}})\cong K^{b}(A\mbox{-}\mathcal{G}proj)$, it
follows that $D^{b}_{gp}(A\mbox{-}{\rm mod})/K^{b}(A\mbox{-}\mathcal{G}proj)=0$. By Lemma 3.1 we get that $D^{b}_{defect}(A)=0$. [BJO, Theorem 2.8] implies
$A$ is Gorenstein. \hfill $\blacksquare$

\vskip 10pt

\begin{cor} \ Let $A$ be a virtually Gorenstein algebra or of finite CM-type. Then  $D^{b}_{gp}(A\mbox{-}{\rm{mod}})$ has Auslander-Reiten triangles if and only if $A$ is Gorenstein.
\end{cor}

\vskip 5pt

\noindent{\bf Proof.}\quad  First we define the functor $\phi: A\mbox{-}\mathcal{G}proj\to {\rm mod}(A\mbox{-}\mathcal{G}proj)\mbox{-}{\rm proj}$ by $\phi(G):={\rm Hom}_{A}(-,G)$ for all $G\in A\mbox{-}\mathcal{G}proj$. Then by Yoneda Lemma we know that $\phi$ is an equivalence. This follows that $K^{b}(A\mbox{-}\mathcal{G}proj)\cong K^{b}({\rm mod}(A\mbox{-}\mathcal{G}proj)\mbox{-}{\rm proj})$ and $K^{-,gpb}(A\mbox{-}\mathcal{G}proj)\cong K^{-,b}({\rm mod}(A\mbox{-}\mathcal{G}proj)\mbox{-}{\rm proj})$.

\vskip 5pt

Whenever $A$ is a virtually Gorenstein algebra or of finite CM-type,  by [Be4] $A\mbox{-}\mathcal{G}proj$ is contravariantly
finite in $A\mbox{-}{\rm mod}$. This follows from [GZ, Theorem 3.6] that  $D^{b}_{gp}(A\mbox{-}{\rm{mod}})\cong K^{-,gpb}(A\mbox{-}\mathcal{G}proj)$. By  Corollary 3.2 we get that $A$ is Gorenstein if and only if there is a triangle-equivalence $ K^{-,b}({\rm mod}(A\mbox{-}\mathcal{G}proj)\mbox{-}{\rm proj})\cong K^{b}({\rm mod}(A\mbox{-}\mathcal{G}proj)\mbox{-}{\rm proj})$. By [Ha3] we get that $K^{-,b}({\rm mod}(A\mbox{-}\mathcal{G}proj)\mbox{-}{\rm proj})$ has Auslander-Reiten triangles  if and only if the category ${\rm mod}(A\mbox{-}\mathcal{G}proj)$ has finite global dimension if and only if $ K^{-,b}({\rm mod}(A\mbox{-}\mathcal{G}proj)\mbox{-}{\rm proj})\cong K^{b}({\rm mod}(A\mbox{-}\mathcal{G}proj)\mbox{-}{\rm proj})$. Hence
 $D^{b}_{gp}(A\mbox{-}{\rm mod})$ has Auslander-Reiten triangles if and only if
$A$ is Gorenstein by above arguments. \hfill $\blacksquare$

\vskip 10pt

Naturally there are two ways to compare Gorenstein derived categories. One way is by Gorenstein derived equivalence. Another way is by recollements of
Gorenstein derived categories. We first show the Gorensteinness is invariant under Gorenstein derived equivalences.

\vskip 10pt

\begin{thm} \ Let $A$ and $B$ be two artin algebras such that $A\mbox{-}\mathcal{G}proj$ and $B\mbox{-}\mathcal{G}proj$ are contravariantly finite respectively. If $A$ and $B$ are  Gorenstein derived equivalent, then  $A$ is Gorenstein if and only if $B$ is Gorenstein.
\end{thm}

\vskip 5pt

\noindent{\bf Proof.} \ Since $A$ and $B$ are  Gorenstein derived equivalent, it follows from [GZ, Proposition 4.2] that
there is a triangle-equivalence $$K^{b}(A\mbox{-}\mathcal{G}proj)\cong K^{b}(B\mbox{-}\mathcal{G}proj).$$
Thus by Lemma 3.1 we get a triangle-equivalence
$$D^{b}_{defect}(A)\cong D^{b}_{defect}(B).$$
This implies that $A$ is Gorenstein if and only if $B$ is Gorenstein. \hfill$\blacksquare$

\vskip 10pt

Let $A$ be an algebra of finite CM-type over a commutative artin ring $R$.  Not much is known about its relative Auslander algebra. Next we will show that the torsionless-finiteness  is invariant under stable equivalences. We first need the following lemmas.

\vskip 10pt

\begin{lem} \
Let $A$ be an artin $R$-algebra of finite CM-type. Then $A\mbox{-}{\mathcal G}proj$ is a dualizing $R$-variety, also ${\rm mod}(A\mbox{-}{\mathcal G}proj)$ is  a dualizing $R$-variety.
\end{lem}

\vskip 5pt

\noindent{\bf Proof.}\quad The fact that $A$ is an artin algebra shows that $A\mbox{-}{\mathcal G}proj$ is a finite $R$-variety. Since $\mathcal{G}p(A)$ is an artin ring, we know that  ${\rm mod}\mbox{-}\mathcal{G}p(A)$ is an abelian category. Hence by isomorphisms
${\rm mod}(A\mbox{-}{\mathcal G}proj)$ $\cong {\rm mod}(\mathcal {G}p(A)\mbox{-}{\rm proj})$ $\cong{\rm mod}\mbox{-}\mathcal{G}p(A)$
we get that ${\rm mod}(A\mbox{-}{\mathcal G}proj)$ is abelian and so $A\mbox{-}{\mathcal G}proj$ has pseudokernels. Next, suppose  $E$ is in $A\mbox{-}{\mathcal G}proj$. Then the object $G$ in ${\rm mod}(A\mbox{-}{\mathcal G}proj)$ obviously has the property that
${\rm Hom}_{A}(X, E)\cong {\rm Hom}_{\mathcal {G}p(A)}({\rm Hom}_{A}(G, $ $X), {\rm Hom}_{A}(G, E))$ for all $X$ in ${\rm mod}(A\mbox{-}{\mathcal G}proj)$.
In order to finish the proof we observe that $(A\mbox{-}{\mathcal G}proj)^{\rm op}$ is equivalent to
$A^{\rm op}\mbox{-}{\mathcal G}proj$ by means of the duality $A\mbox{-}{\mathcal G}proj\to A^{\rm op}\mbox{-}{\mathcal G}proj$ given by $E\mapsto {\rm Hom}_{A}(E, A)$ for all $E$ in $A\mbox{-}{\mathcal G}proj$. Since $A^{\rm op}$ is also an artin algebra with center $R$, it
follows that $(A\mbox{-}{\mathcal G}proj)^{\rm op}=A^{\rm op}\mbox{-}{\mathcal G}proj$ has the properties just derived for $A\mbox{-}{\mathcal G}proj$.
Therefore $A\mbox{-}{\mathcal G}proj$ satisfies the conditions of Theorem 2.4 in [AR] and so
$A\mbox{-}{\mathcal G}proj$ is a dualizing $R$-variety. This implies that ${\rm mod}(A\mbox{-}{\mathcal G}proj)$ is also a dualizing $R$-variety by [AR, Proposition 2.6]. \hfill $\blacksquare$

\vskip 10pt

\begin{lem} \  Let $A$ and $B$ be two artin $R$-algebras of finite CM-type. If $\mathcal{G}p(A)$ and $\mathcal{G}p(B)$ are
stably equivalent, then

\vskip 5pt

$(1)$ \ there is a bijection between the isomorphism classes of indecomposable non-injective torsionless modules in ${\rm mod}\mbox{-}\mathcal{G}p(A)$ and those
in ${\rm mod}\mbox{-}\mathcal{G}p(B)$.

\vskip 10pt

$(2)$ \ there is a bijection between the isomorphism classes of indecomposable nonsimple non-projective injective modules in ${\rm mod}\mbox{-}\mathcal{G}p(A)$ and those in ${\rm mod}\mbox{-}\mathcal{G}p(B)$.
\end{lem}

\vskip 5pt

\noindent{\bf Proof.}\quad   Since $A$ and $B$ are of finite CM-type, it follows from Lemma 3.4 that $A\mbox{-}{\mathcal G}proj$ and $B\mbox{-}{\mathcal G}proj$ are  dualizing $R$-varieties.  Since $\mathcal{G}p(A)$ and $\mathcal{G}p(B)$ are
stably equivalent, we have that $A\mbox{-}{\mathcal G}proj$ and $B\mbox{-}{\mathcal G}proj$ are two stably equivalent dualizing
R-varieties. The result follows from [AR, Corollary 9.11 and 9.14].

\vskip 10pt

Following Ringel [Ri], an artin algebra $A$ is torsionless-finite provided there are only finitely many isomorphism classes of indecomposable torsionless $A$-modules.

\vskip 10pt

\begin{thm} \ Let $A$ and $B$ be two  artin $R$-algebras of finite CM-type. If $A\mbox{-}{\rm mod}/A\mbox{-}{\mathcal G}proj$ and  $B\mbox{-}{\rm mod}/B\mbox{-}{\mathcal G}proj$ are equivalent as categories, then  $\mathcal{G}p(A)$ is torsionless-finite if and only if $\mathcal{G}p(B)$ is torsionless-finite.
\end{thm}

\vskip 5pt

\noindent{\bf Proof.}\quad  By [G1] we know that if $A\mbox{-}{\rm mod}/A\mbox{-}{\mathcal G}proj$ and  $B\mbox{-}{\rm mod}/B\mbox{-}{\mathcal G}proj$ are equivalent as categories, then $\mathcal{G}p(A)$ and $\mathcal{G}p(B)$ are stably equivalent. By Lemma 3.6 we immediately deduce that $\mathcal{G}p(A)$ is torsionless-finite if and only if $\mathcal{G}p(B)$ is torsionless-finite.

\vskip 10pt

Next, we will compare the invariants along recollements of Gorenstein derived categories. We shows the Gorensteinness, the finiteness of Cohen-Macaulay type
and finistic dimension  are invariant with respect to such recollements. We need the following lemma, which is inspired by [AKLY].

\vskip 10pt

\begin{lem} \ Let $A$ and $B$ be two artin algebras and  $F: K^{-,gpb}(A\mbox{-}\mathcal{G}proj)\to K^{-,gpb}(B\mbox{-}$ $\mathcal{G}proj)$   a triangle functor. The following two conditions are equivalent

\vskip 5pt

$(1)$ \ $F$ restricts to $K^{b}(A\mbox{-}\mathcal{G}proj)$;

\vskip 5pt

$(2)$ \ $F(G)\in K^{b}(B\mbox{-}\mathcal{G}proj)$ for any $G\in A\mbox{-}\mathcal{G}proj$.
\end{lem}

\vskip 10pt

\begin{thm} \ Let $A$, $B$ and $C$ be virtually Gorenstein algebras. Assume that  $D^{b}_{gp}(A\mbox{-}{\rm mod})$ admits the following
recollement
$$
D^{b}_{gp}(B\mbox{-}{\rm mod}) \begin{smallmatrix}
  \underleftarrow{ \ \ \  i^{*} \ \ \  } \\
  \underrightarrow{ \ \ \ i_{*} \ \ \ } \\
  \overleftarrow{ \ \ \ i^{!} \ \ \ }
\end{smallmatrix}
D^{b}_{gp}(A\mbox{-}{\rm mod})
\begin{smallmatrix}
  \underleftarrow{ \ \ \ j_{!} \ \ \ } \\
  \underrightarrow{\ \ \ j^{*} \ \ \ } \\
  \overleftarrow{ \ \ \ j_{*} \ \ \ }
\end{smallmatrix}
D^{b}_{gp}(C\mbox{-}{\rm mod})
$$
Then

\vskip 5pt

$(1)$ \ $A$ is of finite CM-type if and only if $B$ and $C$ are so;

\vskip 5pt

$(2)$ \ $A$ is Gorenstein if and only if $B$ and $C$ are so.
\end{thm}

\vskip 5pt

\noindent{\bf Proof.}\quad \  Since $A$ is virtually Gorenstein, it follows from [Be4] that $A\mbox{-}{\mathcal G}proj$ is contravariantly finite in $A\mbox{-}{\rm mod}$. This implies that $D^{b}_{gp}(A)\cong K^{-,gpb}(A\mbox{-}\mathcal{G}proj)$ by the proof of [GZ, Theorem 3.6(ii)]. Similar for $B$ and $C$. By the proof of [GZ, Proposition 4.2 and Lemma 4.3] we get that
$i^{*}: D^{b}_{gp}(A\mbox{-}{\rm mod})\to D^{b}_{gp}(B\mbox{-}{\rm mod})$ and $ j_{!}: D^{b}_{gp}(C\mbox{-}{\rm mod})\to D^{b}_{gp}(A\mbox{-}{\rm mod})$   can restrict to $K^{b}(A\mbox{-}\mathcal{G}proj)$ and $K^{b}(C\mbox{-}\mathcal{G}proj)$ respectively. since
$(i^{*}, i_{*})$ is an adjoint pair and $i^{*}$ is a functor from $D^{b}_{gp}(A\mbox{-}{\rm mod})$ to $D^{b}_{gp}(B\mbox{-}{\rm mod})$, we easily get that $i_{*}$
restricts to $K^{b}(A\mbox{-}\mathcal{G}proj)$.

 \vskip 5pt

If $A$ is of finite CM-type, then by Lemma 3.8 and [Ni] $i^{*}(G_{A})\in K^{b}(B\mbox{-}\mathcal{G}proj)$ is a generator of $D^{b}_{gp}(B\mbox{-}{\rm mod})$. This implies that $B$ is of finite CM-type.  For the converse, by above arguments and Lemma 3.8 we get that $i_{*}(G_{B})$ and $ j_{!}(G_{C})$ are in $K^{b}(A\mbox{-}\mathcal{G}proj)$. Notice from [Ni] that $i_{*}(G_{B})$ and $ j_{!}(G_{C})$ generate $D^{b}_{gp}(A)$. This follows that $A$ is of finite CM-type.

\vskip 5pt

If $A$ is Gorenstein, then  by [GZ, Corollary 3.8] we have that $D^{b}_{gp}(A)\cong K^{b}(A\mbox{-}\mathcal{G}proj)$. Since $i^{*}(G)\in K^{b}(B\mbox{-}\mathcal{G}proj)$ for any $G\in A\mbox{-}\mathcal{G}proj$ by Lemma 3.8, we get that every finitely-generated $B$-module $M$ admits a proper Gorenstein-projective resolution of finite length. This means  $D^{b}_{gp}(B)\cong K^{b}(B\mbox{-}\mathcal{G}proj)$. Hence by Corollary 3.2 we get that $B$ is Gorenstein. For the converse, if $B$ and $C$ are Gorenstein, then we can similarly prove  that $D^{b}_{gp}(A)\cong K^{b}(A\mbox{-}\mathcal{G}proj)$. This implies that $A$ is Gorenstein. \hfill $\blacksquare$

\vskip 10pt

Now we will show the finiteness of finitistic dimension is invariant with respect to the recollements of bounded Gorenstein derived categories.

\vskip 10pt

\begin{lem}\  Let $A$ be an artin of finite CM-type. Then ${\rm fd}(A)<\infty$ if and only if ${\rm
fd}(\mathcal{G}p(A))<\infty$.
\end{lem}

\noindent{\bf Proof.} \ Let ${\rm fd}(A)<\infty$.  Let $_{\mathcal{G}p(A)}M$ be a
$\mathcal{G}p(A)\mbox{-}$module of finite projective dimension and
let $_{\mathcal{G}p(A)}N$ be an arbitrary
$\mathcal{G}p(A)\mbox{-}$module. Let $P^{\bullet}_{M}$ and
$P^{\bullet}_{N}$ be projective resolutions of $M$ and $N$. Note
that by assumption $P^{\bullet}_{M}\in
K^{b}(\mathcal{G}p(A)\mbox{-}{\rm proj})$ and $P^{\bullet}_{N}\in
K^{-,b}(\mathcal{G}p(A)\mbox{-}{\rm proj})$. Then we have
isomorphisms for all $n\in \mathbb{N}$
\begin{align*}
{\rm{Ext}}^{n}_{\mathcal{G}p(A)}(M, N) &=
{\rm{Hom}}_{D^b(\mathcal{G}p(A)\mbox{-}{\rm mod})}(M,
N[n])\\&\cong{\rm{Hom}}_{K^{-,b}(\mathcal{G}p(A)\mbox{-}{\rm
proj})}(P^{\bullet}_{M}, \
P^{\bullet}_{N}[n])\\&\cong{\rm{Hom}}_{K^{-,gpb}(A\mbox{-}\mathcal{G}proj)}(G\otimes_{\mathcal{G}p(A)}P^{\bullet}_{M},
\
G\otimes_{\mathcal{G}p(A)}P^{\bullet}_{N}[n])\\&\cong{\rm{Hom}}_{D^b_{gp}(A\mbox{-}{\rm
mod})}(G\otimes_{\mathcal{G}p(A)}P^{\bullet}_{M}, \
G\otimes_{\mathcal{G}p(A)}P^{\bullet}_{N}[n]).
\end{align*}
Consider the complexes
$E^{\bullet}_{1}=(E^{i}_{1},
d^i)=G\otimes_{\mathcal{G}p(A)}P^{\bullet}_{M}$ and
$E^{\bullet}_{2}=G\otimes_{\mathcal{G}p(A)}P^{\bullet}_{N}$. Note
that $E^{i}_{2}=0$ for $i>0$ and that  ${\rm H}^{-s}{\rm Hom}_{A}(G, E^{\bullet}_{1})=0$
for $s>1$. This follows that ${\rm H}^{-s}{\rm Hom}_{A}(G^{'}, E^{\bullet}_{1})=0$
for $s>1$ and  each  Gorenstein-projective module $G^{'}$. Set
$$X^{\bullet}:=\cdots \to 0\to {\rm Im}d^{-3}\to E_{1}^{-2}\to E_{1}^{-1}\to E_{1}^{0}\to 0\to \cdots.$$
in $D^{b}_{gp}(A\mbox{-}{\rm mod})$ such that $X^{\bullet}\cong E^{\bullet}_{1}$.
Since $X^{\bullet}\cong E^{\bullet}_{1}\in
K^{b}(A\mbox{-}\mathcal{G}proj)$, we infer that $\mathcal{GP}{\rm
dim} _{A}\ {\rm Im}d^{-3}<\infty$. Therefore $X^{\bullet}\cong
E^{\bullet}\in K^{b}(A\mbox{-}\mathcal{G}proj)$ with width
$w(E^{\bullet})\leq {\rm fGd}(A)+4$. Thus $E^{\bullet}$ and
$E^{\bullet}_{2}[n]$ have disjoint support for $n>{\rm fGd}(A)+4$.
The isomorphisms above then shows that
${\rm{Ext}}^{n}_{\mathcal{G}p(A)}(M, N)=0$ for $n>{\rm fGd}(A)+4$,
\ so  ${\rm{Ext}}^{n}_{\mathcal{G}p(A)}(M, N)=0$ for $n>{\rm fd}(A)+4$  by Theorem 2.28 in [Hol].
This implies ${\rm fd}(\mathcal{G}p(A))<\infty$.

\vskip 5pt

Let ${\rm fd}(\mathcal{G}p(A))<\infty$. Let $G^{\bullet}_{M}$ and
$G^{\bullet}_{N}$ be  proper Gorenstein-projective resolutions of
arbitrary $A\mbox{-}$modules $M$ and $N$. Then we have isomorphisms
for $n\in \mathbb{Z}$
\begin{align*}
{\rm Ext}_{\mathcal{G}}^{n}(M,N)&={\rm Ext}_{A\mbox{-}\mathcal{GP}}^{n}(M,N)\\&={\rm
Hom}_{D^b_{gp}(A\mbox{-}{\rm mod})}(M,N[n])\\&\cong {\rm
Hom}_{D^b_{gp}(A\mbox{-}{\rm mod})}(G^{\bullet}_{M}, G^{\bullet}_{N}[n])\\&\cong {\rm
Hom}_{D^b(\mathcal{G}p(A)\mbox{-}{\rm mod})}({\rm Hom}_{A}(G,
G^{\bullet}_{M}),\ {\rm Hom}_{A}(G, G^{\bullet}_{N}[n]))\\&\cong
{\rm Hom}_{D^b(\mathcal{G}p(A)\mbox{-}{\rm mod})}({\rm Hom}_{A}(G,
M),\ {\rm Hom}_{A}(G, N[n]))\\&\cong {\rm
Ext}_{\mathcal{G}p(A)}^{n}({\rm Hom}_{A}(G, M),\
{\rm Hom}_{A}(G, N))
\end{align*}

Since ${\rm fd}(\mathcal{G}p(A))<\infty$, we get that ${\rm
Ext}_{\mathcal{G}}^{n}(M,N)=0$ for all $n>{\rm
fd}(\mathcal{G}p(A))$. This follows that ${\rm fd}(A)={\rm
fGd}(A)<\infty$  again by Theorem 2.28 in [Hol]. \hfill $\blacksquare$

\vskip 10pt

The following corollary is proved by Asadollahi, Hafezi and Vahed in [AHV]. Now we give an alternative proof.

\vskip 10pt

\begin{cor}\ Let $A$ and $B$ be  artin algebras of
finite CM-type. If $A$ and $B$ are Gorenstein derived equivalent,
then ${\rm fd}(A)<\infty$ if and only if ${\rm fd}(B)<\infty$.
\end{cor}

\noindent{\bf Proof.} \   Since $ D^{b}_{gp}(A\mbox{-}{\rm
mod})\cong D^{b}_{gp}(B\mbox{-}{\rm mod})$ by assumption, we get
that $\mathcal{G}p(A)$ and $\mathcal{G}p(B)$ are derived equivalent.
This follows that
 ${\rm fd}(\mathcal{G}p(A))<\infty$ if and only if ${\rm
fd}(\mathcal{G}p(B))<\infty$ by [PX, Theorem 1.1].

\vskip 5pt

By Lemma 3.10 we get that ${\rm fd}(A)<\infty$ if and only if
${\rm fd}(\mathcal{G}p(A))<\infty$, and ${\rm fd}(B)<\infty$ if and
only if ${\rm fd}(\mathcal{G}p(B))<\infty$. Thus ${\rm
fd}(A)<\infty$ if and only if ${\rm fd}(B)<\infty$. \hfill
$\blacksquare$

\vskip 10pt

\begin{thm}\  Let $A$ be an artin algebra of finite CM-type.
 Assume that $D^{b}_{gp}(A\mbox{-}{\rm mod})$ has a
 recollement relative to $D^{b}_{gp}(B\mbox{-}{\rm mod})$ and $D^{b}_{gp}(C\mbox{-}{\rm mod})$
 for artin algebras $B, C$ of finite CM-type.
 Then ${\rm fd}A<\infty$ if and only if
\  ${\rm fd}B<\infty$ and ${\rm fd}C<\infty$.
\end{thm}

\noindent{\bf Proof.} \  Since $A, \ B$ and $C$ are
artin algebras of finite CM-type, there are
triangle-equivalences
 $D^{b}_{gp}(A\mbox{-}{\rm mod})\cong D^{b}(\mathcal{G}p(A)\mbox{-}{\rm
 mod})$, \ $ D^{b}_{gp}(B\mbox{-}{\rm mod})\cong D^{b}(\mathcal{G}p(B)\mbox{-}{\rm
 mod})$  and also $ D^{b}_{gp}(C\mbox{-}{\rm mod})$ $\cong D^{b}(\mathcal{G}p(C)\mbox{-}{\rm
 mod})$. It follows that $D^{b}(\mathcal{G}p(A)\mbox{-}{\rm
 mod})$ has a recollement relative to $D^{b}(\mathcal{G}p(B)\mbox{-}{\rm
 mod})$ and $D^{b}(\mathcal{G}p(C)\mbox{-}{\rm
 mod})$. Hence  by Theorem 3.3 in
[Ha4] ${\rm fd}(\mathcal{G}p(A))<\infty$ if and only if  ${\rm
fd}(\mathcal{G}p(B))<\infty$ and ${\rm
fd}(\mathcal{G}p(C))<\infty$.

\vskip 5pt

By Lemma 3.10  we get that ${\rm fd}A<\infty$ if and only if
${\rm fd}(\mathcal{G}p(A))<\infty$, \  ${\rm fd}B<\infty$ if and
only if ${\rm fd}(\mathcal{G}p(B))<\infty$,   and  ${\rm
fd}C<\infty$ if and only if ${\rm
fd}(\mathcal{G}p(C))<\infty$, respectively. This completes the
proof. \hfill $\blacksquare$

\vskip 10pt

\section {\bf On Gorenstein defect categories}

\vskip 10pt

The previous section used Gorenstein defect categories as a crucial tool. However, not much is known about Gorenstein defect categories. In this section, we will show the Karoubianness of Gorenstein defect categories, and establish relations between  Gorenstein defect categories and Gorenstein stable categories.

\vskip 10pt

We first determine the dimension of  Gorenstein defect categories for a simple class of algebras.

\vskip 10pt

\begin{exm} \ Let $A$ be a representation-finite artin algebra. Then ${\rm dim}D^{b}_{defect}(A)\leq 1$.
\end{exm}

\vskip 5pt
\noindent{\bf Proof.}\quad Since $A$ is representation-finite, it follows from [O] that ${\rm dim}D^{b}(A)\leq 1$.
By [Ro, Lemma 3.4] we get that ${\rm dim}D^{b}_{defect}(A)\leq {\rm dim}D_{sg}(A)\leq {\rm dim}D^{b}(A)$. Hence
${\rm dim}D^{b}_{defect}(A)$ $\leq 1$.
\hfill $\blacksquare$

\vskip 10pt

Now we study the Karoubianness  of the  Gorenstein defect category of an algebra. We need some preparation.

\vskip 10pt

Suppose $A$ is of finite CM-type. For any $M$ in  $A\mbox{-}{\rm mod}/A\mbox{-}\mathcal{G}proj$, we can take a right $A\mbox{-}{\mathcal G}proj\mbox{-}$approximation of $M$, and denote
its kernel by $M^{1}$. Then we have a functor $\overline{\Omega}: A\mbox{-}{\rm mod}/A\mbox{-}\mathcal{G}proj$ $\to A\mbox{-}{\rm mod}/A\mbox{-}\mathcal{G}proj$, view $\overline{\Omega}(M):=M^{1}$.  Denote by $\overline{\Omega}^{n}$ the $n$th composition functor of $\overline{\Omega}$
for any positive integer $n\geq 2$. Then we have

\vskip 10pt

\begin{lem} \  Let $X^{\bullet}$ be a complex in $D^{b}_{gp}(A\mbox{-}{\rm mod})/K^{b}(A\mbox{-}\mathcal Gproj)$ and $n_{0}>0$. Then for any
$n$ large enough, there exists a module $M$ in $\overline{\Omega}^{n_0}(A\mbox{-}{\rm mod})$ such that
$X^{\bullet}\simeq Q(M)[n]$, where $Q: D_{gp}^{b}(A\mbox{-}{\rm mod})\to D^{b}_{gp}(A\mbox{-}{\rm mod})/K^{b}(A\mbox{-}\mathcal Gproj)$ is the quotient functor.
\end{lem}

\vskip 5pt

\noindent{\bf Proof.}\quad Take  an $A\mbox{-}\mathcal Gproj$-quasi-isomorphism  $G^{\bullet}\to X^{\bullet}$
with $G^{\bullet}\in K^{-}(A\mbox{-}\mathcal Gproj)$. Take $n\geq n_{0}$
such that $H^i{\rm Hom}_A(E, X^{\bullet})=0$ for all $i\leq n_{0}-n$ and $E\in A\mbox{-}\mathcal{G}proj$.
Consider the truncation $\tau^{\geq -n}G^{\bullet}=\cdots \to 0 \to M \to G^{1-n}\to G^{2-n}\to \cdots$
of $G^{\bullet}$, which is $A\mbox{-}\mathcal Gproj$-quasi-isomorphic to $G^{\bullet}$. Then the cone of the obvious
chain map $\tau ^{\geq -n}G^{\bullet}\to M[n]$ is in $K^b(A\mbox{-}\mathcal{G}proj)$, which becomes an
isomorphism in $D^{b}_{gp}(A\mbox{-}{\rm mod})/K^{b}(A\mbox{-}\mathcal Gproj)$. This shows that $X^{\bullet}\cong Q(M)[n]$. We observe that
$M$ lies in $\overline{\Omega}^{n_0}(A\mbox{-}{\rm mod})$.
\hfill $\blacksquare$

\vskip 10pt

\begin{lem} \ Let $0\to M \to G^{1-n}\to \cdots \to G^{0}\to N \to 0$ be an {\it $A\mbox{-}\mathcal Gproj$-acyclic} complex with
each $G^i$ Gorenstein-projective. Then we have an isomorphism $Q(N)\cong Q(M)[n]$ in $D^{b}_{gp}(A\mbox{-}{\rm mod})/K^{b}(A\mbox{-}\mathcal Gproj)$. In particular,
for an $A$-module $M$, we have a natural isomorphism $Q(\overline{\Omega}^{n}(M))\cong Q(M)[-n]$.
\end{lem}

\vskip 5pt

\noindent{\bf Proof.}\quad The stalk complex $N$ is $A\mbox{-}\mathcal Gproj$-quasi-isomorphic to
$\cdots \to 0 \to M \to G^{1-n}\to \cdots \to G^{0}\to 0$. This gives rise to a morphism $N\to M[n]$ in
$D^{b}_{gp}(A\mbox{-}{\rm mod})$, whose cone is in $K^{b}(A\mbox{-}\mathcal{G}proj)$. Then this morphism becomes an isomorphism in $D^{b}_{gp}(A\mbox{-}{\rm mod})/K^{b}(A\mbox{-}\mathcal Gproj)$. \hfill $\blacksquare$

\vskip 10pt

We consider the composite functor $Q': A\mbox{-}{\rm mod}\hookrightarrow D_{gp}^{b}(A$-${\rm mod})\xrightarrow{Q}D^{b}_{gp}(A\mbox{-}{\rm mod})/$ $K^{b}(A\mbox{-}\mathcal Gproj)$. It vanishes on $A\mbox{-}\mathcal{G}proj$, so it induces uniquely a functor
$A\mbox{-}{\rm mod}/A\mbox{-}\mathcal{G}proj \to D^{b}_{gp}(A\mbox{-}{\rm mod})/K^{b}(A\mbox{-}\mathcal Gproj)$, which still denote by $Q'$.
Then   for any modules $M,N$ in $A\mbox{-}{\rm{mod}}/A\mbox{-}\mathcal{G}proj$, the functor $Q'$ induces a natural map by Lemma 4.3
$$\Phi^{0}: {\rm Hom}_A(M,N)/A\mbox{-}\mathcal{G}proj(M,N)\to {\rm Hom}_{D^{b}_{gp}(A\mbox{-}{\rm mod})/K^{b}(A\mbox{-}\mathcal Gproj)}(Q(M),Q(N)).$$
This induces  the map
$$\Phi^{n}: {\rm Hom}_A(\overline{\Omega}^n(M),\overline{\Omega}^n(N))/A\mbox{-}\mathcal{G}proj(M,N)\to {\rm Hom}_{D^{b}_{gp}(A\mbox{-}{\rm mod})/K^{b}(A\mbox{-}\mathcal Gproj)}(Q(M),Q(N))$$
given by $\Phi^{n}(f)=(\theta_N)^{-1}\circ (\Phi^{0}(f)[n])\circ\theta_M$.

\vskip 10pt

Consider the chain of maps ${\rm Hom}_A(\overline{\Omega}^n(M),\overline{\Omega}^n(N))/A\mbox{-}\mathcal{G}proj(M,N)\to
{\rm Hom}_A(\overline{\Omega}^{n+1}(M)$,
$\overline{\Omega}^{n+1}(N))/A\mbox{-}\mathcal{G}proj(M,N)$
induced by $\overline{\Omega}$. Then we have an induced map
$$\Phi:\varinjlim{\rm Hom}_A((\overline{\Omega}^n(M),\overline{\Omega}^n(N))
/A\mbox{-}\mathcal{G}proj(M,N)\to {\rm Hom}(Q(M),Q(N)).$$
\hfill $\blacksquare$

\vskip 10pt

\begin{lem} \  Let $A$ be an artin algebra of finite ${\rm CM}$-type, and $M,N$ be in $A\mbox{-}{\rm{mod}}/A\mbox{-}\mathcal{G}proj$.
Then the map $\Phi$ is an isomorphism.
\end{lem}

\vskip 5pt

\noindent{\bf Proof.}\quad We refer to  [Be1, Theorem 3.8]  for a
detailed proof. \hfill $\blacksquare$

\vskip 10pt

Recall that an additive category $\mathcal {A}$ is Karoubian (i.e. idempotent split) provided that each idempotent
$e: X\to X$ splits, that is, it admits a factorization $X\xrightarrow{u}Y\xrightarrow{v}X$ with
$u\circ v={\rm Id}_Y$. In particular, for an artin algebra $A$, the quotient category
$A\mbox{-}{\rm mod}/A\mbox{-}\mathcal{G}proj$ is Karoubian.

\vskip 10pt

\begin{thm} \ The Gorenstein defect category $D^{b}_{defect}(A)$ of an artin algebra $A$ of finite ${\rm CM}$-type is
Karoubian.
\end{thm}

\vskip 5pt

\noindent{\bf Proof.}\quad  We claim that $D^{b}_{gp}(A\mbox{-}{\rm mod})/K^{b}(A\mbox{-}\mathcal Gproj)$ is Karoubian. By Lemma 4.3 it suffices to show that for each module $M$ in $A\mbox{-}{\rm{mod}}/A\mbox{-}\mathcal{G}proj$, an idempotent
$e: Q(M)\to Q(M)$ splits.  Lemma 4.4 implies that for a large $n$, there is an idempotent
$e^{n}: \overline{\Omega}^n(M)\to \overline{\Omega}^n(M)$ in $A\mbox{-}{\rm mod}/A\mbox{-}\mathcal{G}proj$ which is
mapped by $\Phi$ to $e$. Note that the idempotent $e^n$ splits. Then the idempotent $e$ splits.
By Lemma 3.1 we have a triangle-equivalence $D^{b}_{defect}(A)\cong D^{b}_{gp}(A\mbox{-}{\rm mod})/K^{b}(A\mbox{-}\mathcal Gproj)$.
Hence $D^{b}_{defect}(A)$ is Karoubian.
\hfill $\blacksquare$

\vskip 10pt

Next we will show the Gorenstein defect category of an algebra $A$ is triangular equivalent to
the stabilization of the Gorenstein stable category of $A$.

\vskip 10pt

\begin{lem} \ Let $A$ be an artin algebra such that $A\mbox{-}\mathcal{G}proj$  is contravariantly finite in $A\mbox{-}{\rm mod}$. Then there is a triangle-equivalence
$$S(A\mbox{-}{\rm mod}/A\mbox{-}\mathcal{G}proj)\cong D^{b}_{defect}(A).$$
\end{lem}
\noindent{\bf Proof.} \ Since $A\mbox{-}\mathcal{G}proj$  is contravariantly finite  in
$A\mbox{-}{\rm mod}$,  it follows from [Be1, Theorem 3.8] that  there is a triangle-equivalence
$$S(A\mbox{-}{\rm mod}/A\mbox{-}\mathcal{G}proj)\cong K^{-,gpb}(A\mbox{-}\mathcal{G}proj)/K^{b}(A\mbox{-}\mathcal{G}proj).$$
By [KZ, Finial Remak] we get  a triangle-equivalence
$$S(A\mbox{-}{\rm mod}/A\mbox{-}\mathcal{G}proj)\cong D^{b}_{defect}(A).$$ \hfill$\blacksquare$

\vskip 10pt

As an application we show the equivalences of Gorenstein stable categories can induce the equivalences of Gorenstein defect categories for
two algebras of finite CM-type. For convenience we introduce two definitions.

\vskip 10pt

\begin{defn} \ Two artin algebras $A$ and $B$ are said to be Gorenstein stably equivalent if their Gorenstein stable categories $A\mbox{-}{\rm{mod}}/A\mbox{-}\mathcal{G}proj$ and $B\mbox{-}{\rm{mod}}/B\mbox{-}\mathcal{G}proj$
are equivalent as left triangulated categories.
\end{defn}

\vskip 10pt

\begin{defn} \  Two  artin algebras $A$ and $B$ are said to be Gorenstein defect equivalent if there is a triangle-equivalence
$D^{b}_{defect}(A)\cong D^{b}_{defect}(B)$.
\end{defn}

\vskip 10pt

\begin{cor} \ Let $A$ and $B$ be two artin algebras such that $A\mbox{-}\mathcal{G}proj$ and $B\mbox{-}\mathcal{G}proj$ are contravariantly finite respectively.
 If $A$ and $B$ are Gorenstein stable equivalent, then $A$ and $B$ are Gorenstein defect equivalent. Thus the Gorenstein defect category of
an algebra of finite CM-type is uniquely determined by its Gorenstein stable category.
\end{cor}

\vskip 5pt

\noindent{\bf Proof.} \  Since $A$ and $B$ are Gorenstein stably equivalent, then there is a triangle-equivalence by [Be1, Corollary 3.3]
$$S(A\mbox{-}{\rm mod}/A\mbox{-}\mathcal{G}proj)\cong S(B\mbox{-}{\rm mod}/B\mbox{-}\mathcal{G}proj).$$
By Lemma 4.6 we have  triangle-equivalences
$$S(A\mbox{-}{\rm mod}/A\mbox{-}\mathcal{G}proj)\cong D^{b}_{defect}(A),$$
and
$$ S(B\mbox{-}{\rm mod}/B\mbox{-}\mathcal{G}proj)\cong D^{b}_{defect}(B).$$
Hence we get  a triangle-equivalence
$$D^{b}_{defect}(A)\cong D^{b}_{defect}(B).$$
\hfill$\blacksquare$

\vskip 20pt

{\bf Acknowledgements.} Some of this work was done when the author is visiting Professor Steffen K\"{o}nig in the University of Stuttgart, Germany.
The author would like to thank  Steffen K\"{o}nig  for useful discussions and comments related to this work. The author also would like to thank Professor Pu Zhang for helpful discussions.

\end{document}